\documentclass{article} 
\usepackage{amsmath}
\usepackage{amsfonts}
\usepackage{amsthm}
\usepackage{dsfont}
\usepackage{epsfig}
\usepackage{upgreek}
\usepackage{graphicx}
\usepackage{caption}
\usepackage{setspace}

\theoremstyle{definition}

\newcommand{\Gmi}{\mbox{${\bf G}^{\mbox{\tiny -1}}$}}
\newcommand{\RGA}{\mbox{${\mathcal RGA}$}}

\newcommand{\Amn}{\mbox{${\bf A}$}}
\newcommand{\Amni}{\mbox{${\bf A}^{\!\mbox{\tiny -1}}$}}
\newcommand{\Amnt}{\mbox{${\bf A}^{\!\mbox{\tiny T}}$}}
\newcommand{\Amit}{\mbox{${\bf A}^{\!\mbox{\tiny -T}}$}}
\newcommand{\Bm}{\mbox{${\bf B}$}}

\newcommand{\sigmin}{\mbox{$\sigma_{\mbox{\tiny min}}$}}

\newcommand{\Rm}{\mbox{${\bf R}$}}

\newcommand{\Rmt}{\mbox{${\bf R}^{\mbox{\tiny T}}$}}
\newcommand{\Vm}{\mbox{${\bf V}$}}
\newcommand{\Vmt}{\mbox{$\transp{\Vm}$}}
\newcommand{\Umt}{\mbox{$\transp{\Um}$}}

\newcommand{\Mm}{\mbox{${\bf M}$}}
\newcommand{\Dm}{\mbox{${\bf D}$}}

\newcommand{\Dmii}{\mbox{${\bf D}^2_{ii}$}}

\newcommand{\Em}{\mbox{${\bf E}$}}
\newcommand{\Emi}{\mbox{${\bf E}^{\mbox{\tiny -1}}$}}
\newcommand{\Emii}{\mbox{${\bf E}^{\mbox{\tiny -1}}_{ii}$}}

\newcommand{\Hm}{\mbox{${\bf H}$}}
\newcommand{\Fm}{\mbox{${\bf F}$}}

\newcommand{\Gm}{\mbox{${\bf G}$}}

\newcommand{\transp}[1]{{#1}^{\mbox{\tiny T}}}
\newcommand{\transpi}[1]{{#1}^{\mbox{\tiny -T}}}

\newcommand{\transpinv}[1]{{#1}^{\mbox{\tiny -T}}}
\newcommand{\inv}[1]{{#1}^{\mbox{\tiny -1}}}

\newcommand{\Um}{\mbox{${\bf U}$}}

\begin{document}

\title{{ Expression of a Real Matrix as a Difference\\
         of a Matrix and its Transpose Inverse}}      
\author{{ Mil Mascaras}\footnote{Work initiated while honorary visiting professor.\newline\newline
Paper submitted to {\em Journal de Ciencia e Ingenieria} on 28 Dec 2018. Accepted with revisions on 2 Apr 2019.
To be published: Vol.\ 11, No.\ 1,  2019.}~ and {Jeffrey Uhlmann}\\
Dept. of Electrical Engineering and Computer Science\\
University of Missouri-Columbia}

\date{}          
\maketitle
\thispagestyle{empty}


\begin{abstract}
In this paper we derive a representation of an arbitrary real 
matrix $\Mm$ as the difference of a real matrix $\Amn$ 
and the transpose of its inverse. This expression may prove 
useful for progressing beyond known results for which the appearance 
of transpose-inverse terms prove to be obstacles, particularly in 
control theory and related applications such as computational
simulation and analysis of matrix representations of
articulated figures.\\
~\\
\begin{footnotesize}
\noindent {\bf Keywords}: {\sf\scriptsize Articulated figure analysis,
algebraic Riccati equation, control systems, 
matrix analysis, matrix decompositions, matrix splitting, octonions, quaternions,
real matrices, relative gain array, RGA, 
singular value decomposition, SVD.}\vspace{4pt}\\ 
\noindent {\scriptsize\sf {\bf AMS Code}: 15A09}
\end{footnotesize}
\end{abstract}

\section{Introduction}

A common way to glean information about a given matrix, e.g.,
to reveal opportunities for manipulating or transforming it, is to
express it as the sum or difference of matrices with particular
structure or properties. This may be a simple expression of a
singular matrix as a sum of nonsingular ones\!~\cite{mathmag} or
a given matrix expressed as the sum of a symmetric matrix and a 
skew-symmetric matrix\!~\cite{HJ1}. {\em Matrix splitting}
is an example of a widely-used technique that relies on the
expression of a matrix as the difference of two matrices with
special properties, e.g., for solving systems of differential 
equations \cite{Varga}. 

Transpose-inverse terms $\transp{(\Amni)}$ 
$\left(\mbox{or, equivalently,\,} \inv{(\Amnt)}\right)$, commonly abbreviated
as $\Amit$, arise naturally in a variety of 
control system contexts, e.g., the relative gain 
array (RGA) \cite{bristol} and formulations of the
controllability Gramian \cite{nature}. As an example,
classical solution methods for the discrete-time
algebraic Riccati equation involve the 
symplectic form \cite{laub}:
\begin{equation}
\left[ \begin{array}{c c}
            \Fm+\Gm\transpi{\Fm}\Hm & -\Gm\transpi{\Fm}\\
            -\transpi{\Fm}\Hm & \transpi{\Fm}
       \end{array}
\right]
\end{equation}
Unfortunately, there are very few matrix factorizations or 
decompositions involving transpose-inverse terms to
aid in the manipulation of equations for theoretical
analysis or practical implementation. In this paper we provide an
incremental improvement to this state of affairs with a representation
of an arbitrary nonsingular real matrix $\Mm$ as the difference
of a real matrix $\Amn$ and its transpose inverse:
\begin{equation}
   \Mm ~=~ \Amn - \Amit .\label{main}
\end{equation}
We begin with a derivation of this result in the next section and then
develop related results for singular and complex matrices. 

\section{The Real Nonsingular Case}

The main result can be derived from an
application of the singular value decomposition (SVD) of the real 
nonsingular matrix $\Mm$ as 
\begin{equation}
      \Mm ~=~ \Um\Dm\Vmt
\end{equation}
where $\Um$ and $\Vm$ are real orthonormal and $\Dm$ is a 
positive diagonal matrix of the singular values of $\Mm$.
It can be observed that defining $\Amn$ with diagonal $\Em$ as 
\begin{equation}
       \Amn ~=~  \Um\Em\Vmt 
\end{equation}
satisfies
\begin{eqnarray}
   \Mm & = & \Amn - \Amit \\
       ~ & = & \Um\Em\Vmt - \transpinv{(\Um\Em\Vmt)}\\ 
       ~ & = & \Um\Em\Vmt - \transp{\left(\inv{(\Um\Em\Vmt)}\right)}\\ 
       ~ & = & \Um\Em\Vmt - \transp{(\Vm\Emi\Umt)}\\ 
       ~ & = & \Um\Em\Vmt - \Um\Emi\Vmt 
\end{eqnarray}
only if each diagonal element $\Em_{ii}$ satisfies
\begin{equation}
    \Dm_{ii} ~=~ \Em_{ii} - \Emii.
\end{equation}
This defines a quadratic constraint on each $\Em_{ii}$
that can be verified to admit solutions
\begin{equation}
      \Em_{ii} ~=~ 
      \frac{1}{2}\left(\Dm_{ii}\pm\sqrt{\Dmii + 4}\right)
        \label{diffconstraints}
\end{equation}
which can be taken as real for all $\Dm_{ii}$ by
positivity of the singular values of $\Mm$. Furthermore,
the nonnegative solution can be taken for each 
$\Em_{ii}$ so that $\Amn$ is determined
completely up to the uniqueness (or lack thereof)
provided by the SVD evaluation method.

The case of singular $\Mm$ can be handled in two ways. 
As it stands, the previous construction can be restricted to 
only the nonzero singular values to give 
\begin{equation}
   \Mm ~=~ \Amn - \transp{(\Amn^{\dagger})} \label{singular}
\end{equation}
where $\Amn$ now has the same rank as $\Mm$ and
its inverse is replaced with a pseudoinverse. Alternatively,
the zero singular values of $\Amn$ above can be replaced with
unity\footnote{This works because the unit
singular values will cancel
in the difference with their corresponding
ones in the transpose inverse, and in fact $1$ is
the positive solution to Eq.\,(\ref{diffconstraints})
for $\Dm_{ii}=0$.} so that 
singular $\Mm$ is expressed as the difference of
nonsingular matrices in the form of Eq.\,(\ref{main}). 
In summary, the expression of Eq.\,(\ref{main})
can be obtained for all real square matrices $\Mm$ whereas
that of Eq.\,(\ref{singular}) extends to real $\Mm$
of any shape.

\section{Variations and Generalizations}

The results from the previous section can be verified
to generalize directly
to the case of complex $\Mm$ if conjugate-transpose is
used in place of the transpose operator
\begin{equation}
   \Mm ~=~ \Amn - (\Amni)^*. \label{complex}
\end{equation}
However, variants in which
the matrix difference of Eq.\,(\ref{main}) is replaced with
a sum, or the transpose operator is maintained for
complex $\Mm$, cannot generally be obtained in a 
similar form.

The challenge posed by the form
\begin{equation}
   \Mm ~=~ \Amn + \Amit \label{sum}
\end{equation}
is that the resulting analog of Eq.\,(\ref{diffconstraints}) 
becomes
\begin{equation}
      \Em_{ii} ~=~ 
      \frac{1}{2}\left(\Dm_{ii}\pm\sqrt{\Dmii-4}\right)
      \label{diffconstraints2}
\end{equation}
which admits real solutions only for $\Dm_{ii}\geq 2$. In other words,
the form of Eq.\,(\ref{sum}) can only be obtained with real $\Amn$
if the smallest singular value of $\Mm$ satisfies $\sigmin\geq 2$.
On the other hand, from this we can conclude that the form of Eq.\,(\ref{sum}) 
can always be obtained for $\frac{2}{\sigmin}\Mm$, i.e., when real 
nonsingular $\Mm$ is scaled to ensure
its smallest nonzero singular value is 2. This yields the slightly less
pleasing form
\begin{equation}
   \Mm ~=~ c\left(\Amn + \Amit\right) \label{sumc}
\end{equation}
which can be obtained for all real nonsingular $\Amn$ 
by letting $c=2/\sigmin$. 

For completeness we note that 
Eq.\,(\ref{diffconstraints}) \& Eq.\,(\ref{diffconstraints2})
can be applied with matrix arguments in place of diagonal
elements to obtain, respectively, the non-transpose forms:
\begin{equation}
   \Mm ~=~ \Amn - \Amni
\end{equation}
and
\begin{equation}
   \Mm ~=~ c\left(\Amn + \Amni\right) 
\end{equation}
which are not relevant to the focus of this paper, e.g., because
they are not applicable to rectangular matrices, but
may be of independent interest.

\section{Application to Articulated Figure Analysis}

For purposes of kinematic simulation and analysis of 
articulated figures, e.g., for motion prediction and/or
animation, there are multiple mathematical representations.
Most commonly, a human figure is represented with
segments and joints (see Fig.\,1) with the mobility 
constraints expressed using Euler angles, quaternions
(or double quaternions), and/or exponential 
maps \cite{Usta,Schilling,Grassia}. 

\begin{figure}[h]
\centering
\includegraphics[width=\textwidth]{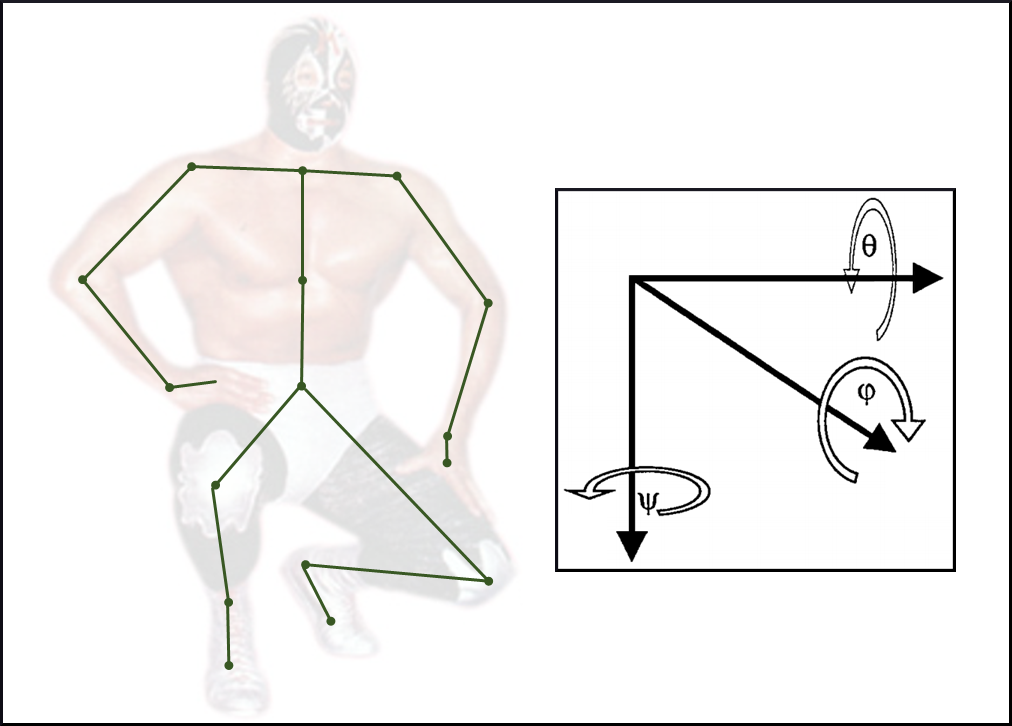}\\
{\em Segment and joint representation of an
articulated human figure.}
\end{figure}

More generally, the structure of the figure may be represented
simply as a matrix. For example, compositions of rigid rotations 
involving a heirarchy of coordinate axes can be represented in 
the form of a real orthonormal rotation matrix, $\Rm$, such that 
the transformation of a given matrix $\Mm$ representing an 
articulated figure can be expressed as $\Rm\Mm\Rmt$.
The price paid for the generality of representing figures
and shapes using matrices is the 
challenge of how to manipulate and analyze such expressions, 
e.g., for operations such as graph matching or simply to
gain enhanced intuitive insights. As indicated in \cite{Mateus}:
``it is not yet clear how to choose and characterize the group of 
transformations under which such shapes should be studied.''

A key property that clearly must be maintained is consistency
with respect to real orthonormal transformations. It can easily
be verified that this property is maintained by the transformation
$f(\Mm)\rightarrow\Amn$, $\Mm=\Amn - \Amit$, by virtue of
its derivation via the SVD, or explicitly as:
\begin{eqnarray}
\Rm\Mm\Rmt & = & \Rm(\Amn - \Amit)\Rmt\\
~      & = & \Rm\Amn\Rmt - \Rm\Amit\Rmt\\
~      & = & (\Rm\Amn\Rmt) - \transpi{(\Rm\Amn\Rmt)}
\end{eqnarray} 
where the fact that $\transpi{\Rm}=\Rm$ for real orthonormal $\Rm$
has been exploited. This therefore demonstrates the desired
consistency property:
\begin{equation}
     f(\Rm\Mm\Rmt) ~=~ \Rm\cdot f(\Mm)\cdot\Rmt.
\end{equation}
In the previous section we discussed the generalization from transpose 
to conjugate-transpose, which can support consistency with respect 
to unitary transformations, and it is natural to consider a further
generalization to the nonassociative octonions. This would permit forces to be more
flexibly incorporated \cite{Weng}, and it would permit temporal
sequences of non-compositional operators to be expressed uniquely
based on a specified associativity rule \cite{Mil}. More specifically,
a temporal sequence of non-associative operators $\alpha_i$ 
applied at times $t_i$ could be derived as a solution to a given
problem and expressed in directional time-assymetric form as
\begin{eqnarray}
   t_0 &\rightarrow& \alpha_0 \nonumber \\
   t_1 &\rightarrow& \alpha_1(\alpha_0)\\  
   t_2 &\rightarrow& \alpha_2(\alpha_1(\alpha_0)) \nonumber\\ 
   ~    &\vdots& ~ \nonumber
\end{eqnarray}
Unfortunately, the proposed matrix decomposition does not appear to
practically accommodate unitary transformations over the octonions 
because there presently does not exist an efficiently computable 
octonion analog of the SVD\footnote{However, the recent approach
of \cite{Mizo} is of potential relevance in this regard.}.

\section{Discussion}

The main result of this paper is that every real
matrix $\Mm$ can be represented as the difference
of two real matrices in the form
\begin{equation}
   \Mm ~=~ \Amn - \Amit . \label{final}
\end{equation}
A measure of the potential utility of a given factorization or 
decomposition is the extent to which it can be used to obtain
nontrivial derivative results. Eq.\,(\ref{final}) admits a variety
of trivial ones, e.g., multiplication of $\Mm$ by $\Amni$ or
$\Amnt$ gives a symmetric difference involving a positive semidefinite
matrix and the identity matrix. However, additional structural
properties can be derived based on known results involving
transpose inverses, e.g., the relative gain array (RGA) 
mentioned in the introduction. The RGA of a square nonsingular
matrix $\Gm$ is defined~\cite{bristol} as
\begin{equation}
   \RGA(\Gm) ~\doteq~ \Gm \circ \transp{(\Gmi)} \label{rga}
\end{equation}
where $\circ$ represents the elementwise Hadamard matrix product.
Because of its practical importance in control applications the RGA 
has been shown~\cite{rgajs} to have a variety of interesting 
mathematical properties, and those properties therefore carry over
to the Hadamard product of the matrix terms of Eq.\,(\ref{final}).
For example, it is known that $\RGA(\Gm)$ is invariant with 
respect to diagonal scalings of $\Gm$ and that its rows and columns
have unit sum, i.e., $\RGA(\Gm)$ is generalized doubly stochastic.
This implies that by letting $\Bm=\Amit$ in Eq.\,(\ref{final}) one
can obtain an expression with significantly less explicit
structural information
\begin{equation}
   \Mm ~=~ \Amn - \Bm,~ ~~ \Amn\circ\Bm\,\in\,\mbox{\em Generalized Doubly Stochastic}
\end{equation}
which would still be of interest based solely on the property
that $\Amn\circ\Bm$ is generalized doubly stochastic\footnote{It has 
recently been shown that a new generalized inverse \cite{simax} can be 
used to preserve some RGA properties in the case of singular 
$\Gm$ \cite{scinv,rga}, thus 
conferring related properties to the expression of
Eq.\,(\ref{singular}).}. In summary, the main result of this paper is a 
specialized form of matrix splitting that offers a diverse
set of interesting mathematical properties that derive from and thus may
prove useful to control theory and related applications.

\end{document}